\documentclass[12pt, english]{article}
\usepackage{amssymb,amsfonts,amsmath,amsthm}
\usepackage{lineno,hyperref}
\usepackage{bbm}
\usepackage{dsfont}
\usepackage{graphicx}
\usepackage{color}
\usepackage{amssymb}
\usepackage{natbib}
\usepackage{xcolor}

\usepackage{float}
\usepackage{color,soul}

\usepackage[margin=1in]{geometry}
\usepackage{enumitem}
\modulolinenumbers[5]

\newcommand{\beqn}{\vspace{-0.25cm}\begin{eqnarray*}}
\newcommand{\eeqn}{\end{eqnarray*}}
\newcommand{\bneqn}{\vspace{-0.25cm}\begin{eqnarray}}
\newcommand{\eneqn}{\end{eqnarray}}

\newcommand{\beq}{\begin{eqnarray*}}
\newcommand{\feq}{\end{eqnarray*}}
\newcommand{\feqn}{\end{eqnarray}}

\newcommand{\expesub}[2]{\mathbb{E}_{#1}\bracks{#2}}
\newcommand{\bracks}[1]{\left[#1\right]}
\newcommand{\parens}[1]{\left(#1\right)}
\newcommand{\expe}[1]{\mathbb{E}\bracks{#1}}

\newtheorem{theorem}{Theorem}
\makeatletter \@addtoreset{theorem}{section}\makeatother

\newtheorem{conjecture}[theorem]{Conjecture}

\newtheorem*{theorem*}{Theorem}

\newtheorem{corollary}[theorem]{Corollary}

\title{Exercising Control When Confronted by a (Brownian) Spider}
\author{Philip Ernst}
\begin{document}
\maketitle

\begin{abstract}
We consider the Brownian ``spider,''  a construct introduced in \cite{Pitman} and in \cite{Dubins}. In this note, the author proves the ``spider'' bounds by using the dynamic programming strategy of guessing the optimal reward function and subsequently establishing its optimality by proving its excessiveness. 
\end{abstract}

\begin{center}
\textit{In memory of my mentor, Professor Larry Shepp (1936-2013)}
\end{center}

\normalsize

\section{Introduction}\label{sec:introduction}
\noindent In this note, we consider the Brownian ``spider,'' a process also known as  the ``Walsh'' Brownian motion, due to \cite{Pitman} and \cite{Dubins}. The Brownian spider is constructed as a set of $n \ge 1$ half-lines, or ``ribs,'' meeting at a common point, $O$. A Brownian motion
on a spider starting at zero may be constructed from a standard reflecting
Brownian motion ($|W_t|,t \geq 0$) by assigning an integer 
$i \in \{1,\ldots,n\}$ uniformly and independently to each excursion
which is then transferred to an excursion on rib $i$ (here, $i$ should be interpreted as the index of the rib on which the next excursion occurs). It is helpful to think about the Brownian spider as a bivariate process; the first coordinate of the process is reflecting Brownian motion and the second coordinate of the process is the rib index. Formally, we define the  Brownian spider process $Z_t$ as
\begin{equation} 
Z_t= \parens{|W_t|, R_t}, t \geq 0
\end{equation}
where $|W_t|$ is reflected Brownian motion and $R_t$ is the rib on which the process is located at time $t$. $|W_t|$ can be decomposed into excursions away from 0 with endpoints $t_k$ s.t. $|W_{t_k}| = 0$. $R_t$ is constant between $t_k$ and $t_{k+1}$ for all $i$, and $R_t = i$ means the excursion occurs on the rib $i$.
We define the supremum of reflected Brownian motion on each rib as
\beqn
S_i(t)= \sup_{\{t:\,\, R_t=i\}}|W_t|, \,\, t \geq 0, \, i=1,...,n.
\eeqn
\indent Below is a sample path realization of the Brownian spider for $n=3$. We use $W_i(t)$ to denote the process on the rib $i$: 
\begin{equation}
W_i(t) = \left\{  \begin{array}{cc}
|W_t|  &   \text{ if }  R_t = i   \\
0   &   \text{ if } R_t \neq i.
\end{array}
\right.
\end{equation}

\begin{figure}[H]
\centering
\includegraphics[width=0.5\linewidth]{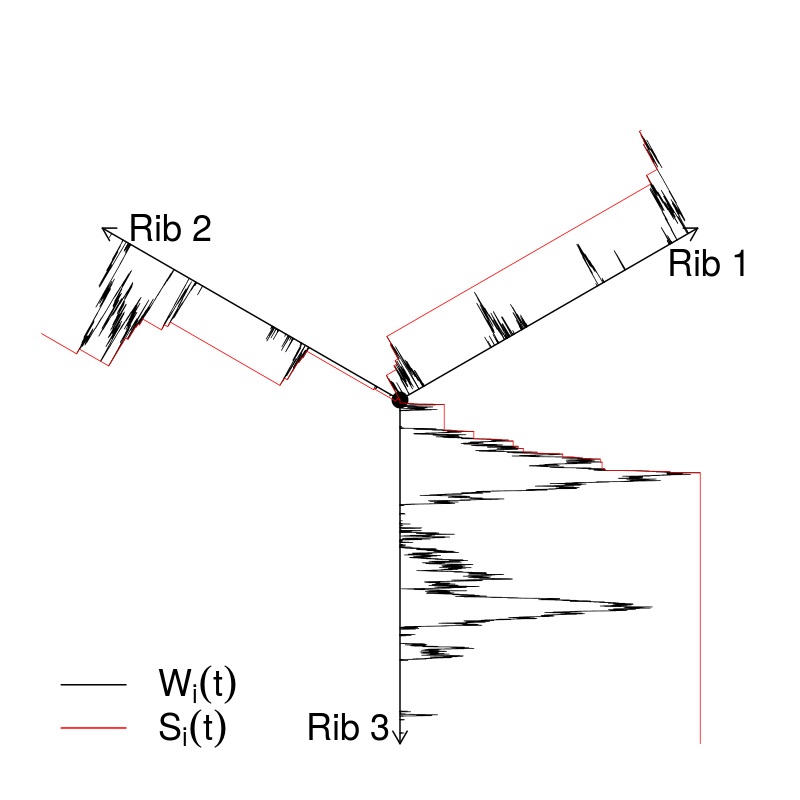}
\caption{A sample path realization of the Brownian spider for $n=3$.}
\label{fig:marginals}
\end{figure}

\indent In an attempt to understand the unboundedness of Brownian motion on the spider up to time $t$, a natural question to ask is: what is $\expe{\sum_{i=1}^nS_i(t)}$?
However, Lester Dubins (personal communication with Larry Shepp) asked a different question. Dubins wished to design a \textit{stopping time} to maximize the coverage of Brownian motion on the spider for a given expected time. That is, he wished to find

\bneqn
C_n := \sup_{\{\tau: \,\, \expe{\tau} =1\}} \expe{S_1(\tau)+\ldots+S_n(\tau)}, 
\eneqn

\noindent
where the supremum is calculated over all stopping times of mean one. Equivalently, Dubins wished to calculate the smallest $C = C_n$ such that for every stopping time $\tau$ the
following inequality holds

\bneqn \label{eq2}
\expe{S_1(\tau)+S_2(\tau)+\ldots+S_n(\tau)} \le C_n \sqrt{\expe{\tau}}.
\eneqn
(note that for any stopping time $\tau$, $\expe{S_i(\tau)}$ scales with with $\sqrt{\tau}$). The left side of equation (\ref{eq2}) is the mean total
measure of space visited on the spider up to time $\tau$. 

\indent For $n = 0$, we,
somewhat inconsistently, define $C_0$ in a similar way for ordinary Brownian
motion without a reflecting barrier at zero. We seek the smallest
constant $C_0$ for which the {\em one-sided maximum} satisfies

\bneqn
\expe{\max_{\{0 \le t \le \tau\}} W(\tau)} \le C_0 \sqrt{\expe{\tau}}.
\eneqn

In this note, we will prove that the optimal bounds $C_n = \sqrt{n+1}$, for $n = 0,1,2$. Without further delay, the author notes that the \textit{solution} of the optimal bounds for $n=0,1,2$ is \textit{not new}. The cases $n=0,1$ were solved by \cite{Dubins}, and the case $n=1$ was also independently solved by \cite{Gilat} by a different method. The $n=2$ case was recently resolved by \cite{Dubins2}. What \textit{is new}, however, is the dynamic programming \textit{strategy} the author employs to find the bounds $C_n$ for $n=0,1,2$, which he believes to be the most tractable approach for solving for $C_n$ for all $n$ (despite much effort by many researchers, this problem remains open). The behavior of $C_n$ for large $n$ is interesting because when $n = \infty$, the total measure of space visited on the spider up to time $t > 0$ is also infinite. This is because it is the total variation of a Brownian motion on $[0,t]$ because at each return to the node, a fresh rib is chosen. 

Larry Shepp saw dynamic programming to be the root of all optimal control problems. In general, there are two strategies that can be used to solve a dynamic programming problem. 

(A) Guess a candidate for an optimal strategy, calculate the reward function for the strategy, then prove its excessiveness.

(B) Guess the optimal reward function and establish its optimality by proving its excessiveness.
 
Unlike \cite{Dubins2}, which employs strategy (A), our approach is that of (B), and to the best of our knowledge, we are the first to do so. In stochastic optimization, strategy (B) reduces to  ``guessing'' the right optimal control function. If one can guess the right function, the supermartingale becomes a martingale, and It\^o calculus  takes care of the rest. This approach appears prominently throughout Shepp's most seminal works, specifically on  p.634 of \cite{Shepp2}, on p.207 of \cite{Shepp3}, p.1528 of \cite{Shepp4}, on p.335 of \cite{Ernst}, and most recently, on p.422 of \cite{Ernst2}.

The organization of this note is as follows: In Section 2, we formalize our guess for the optimal reward function. In Section 3, we establish the optimality of this function by proving its excessiveness, albeit only in the cases $n=0,1,2$. We conclude by arguing the viability of our strategy towards a solution of the general problem.

\section{Our guess of the optimal reward function}\label{sec:formal_model_approach}
Let $r = R_0$ be the index of the starting rib, $x$ be a fixed distance along the rib $r$, and $C$ and $M$ finite constants. In order to obtain the least upper bound $C_n$, we must solve
a more general optimal stopping problem. Let $s_1, s_2,...,s_n \geq 0$ be the distances that have already been covered on each of the respective ribs at time $0$. For {\it every value} of
$C > 0$, and \textit{every} choice of $r$ and $x$ such that $x \le s_r$, and $s_1,\ldots,s_n$, we must find the value of

\bneqn \label{quantity}
V(x,r;s_1,\ldots,s_n;C) := \sup_{\{\tau: \,\, \expe{\tau} \leq M\}} \expesub{\{x,r,s_1,\ldots,s_n\}}{S_1(\tau)+\ldots + S_n(\tau) - C\tau}.
\eneqn

The subscript of the expectation, ${\{x,r,s_1,\ldots,s_n\}}$, denotes that the process is currently at a distance $x$ on rib $r$ at time $0$. By abuse of notation, $S_i(\tau)$ denotes the furthest point covered on rib $i$ up to time $\tau$. Note that we must find $V(x,r;s_1,\ldots,s_n;C)$ not only for $x=0$ and $s_1 = \ldots = s_n = 0$, but for
{\it every point} of the spider at $x$ on {\it every} rib $r$ as initial
point, and every starting position for $s_i, i=1,\ldots,n$.

In (\ref{quantity}), the supremum is taken over bounded stopping times $\tau$. Even though we
only need the case when the initial point is $O$ and when $s_i = 0,\,\,i=1,\ldots,n$, standard
martingale methods of solving optimal stopping problems do not work unless we
can find the formula $V$ for every starting position (see, for example: ~\cite{Brumelle},~\cite{Dynkin},  ~\cite{Shepp}, ~\cite{Talagrand}, and ~\cite{Walrand}). 

We ``guess''  that ${\hat{V}}(x,r,s_1,\ldots,s_n,C)$ should have the following properties:
\medskip

(a) ${\hat{V}}(0,r,s_1,\ldots,s_n,C)$ does not depend on $r$ (if $x_i=0, r$ becomes irrelevant).
\medskip

(b) $\frac{d}{dx}{\hat{V}}(x,r,s_1,\ldots,s_n,C) = 0$ at $x=0$ $\forall r$.
\medskip

(c) $\frac{d}{ds_r}{\hat{V}}(x,r,s_1,\ldots,s_n,C) = 0$ at $x = s_r\,\, \forall \,r$.
\medskip

(d) $\frac{1}{2} \frac{d^2}{dx^2}{\hat{V}}(x,r,s_1,\ldots,s_n,C) \leq C, \, 0 \le x \le s_r,\,r=1,\ldots,n$.
\medskip

(e) ${\hat{V}}(x,r,s_1,\ldots,s_n,C) \ge s_1+\ldots+s_n$, for all $0 \le x \le s_r$ and $r = 1,\ldots,n$.
\medskip

(f) If strict inequality holds in property (e), $\frac{1}{2} \frac{d^2}{dx^2}{\hat{V}}(x,r,s_1,\ldots,s_n,C) = -C$, $0 \le x \le s_r$.\\

\noindent Intuitively, at a stopping place, we are far from any boundary point where an
$s$ would increase and thus we are willing to accept the reward
${\hat{V}}(x,r,s_1,\ldots,s_n,C) = s_1+\ldots+s_n$. 

\section{Establishing the optimality of the reward function}\label{sec:formal_model_approach}

\begin{theorem}
If we have a function ${\hat{V}}$ satisfying properties (a)-(f) in Section 2, then

\bneqn
V(x,r,s_1,\dots,s_n;C) \equiv {\hat{V}}(x,r,s_1,\dots,s_n,C).
\eneqn
\end{theorem}
\begin{proof} 
Consider the process

\bneqn
Y(t) = {\hat{V}}\parens{Z_t,{\bf S}(t),C} - Ct, \,\,\,t \ge 0
\eneqn
where $\mathbf{S}(t) = \parens{S_1(t),...,S_n(t)}$. $Y(t)$ is a continuous local supermartingale at
$x=0$ by properties (a) and (b), at $x = s_r$ by property (c), and at any $x$ by property (d). For any bounded stopping time $\tau$, it follows from the optional sampling theorem that $\expe{Y(\tau)} \leq Y(0)$.
Property (e) gives us that for any bounded $\tau$, 
\bneqn \label{coolone}
\expesub{\{x,r,s_1,\ldots,s_n\}}{S_1(\tau)+\ldots+S_n(\tau)- C \tau} \le {\hat{V}}\parens{x,r,s_1,\ldots,s_n,C}.
\eneqn
From the definition of $V$ in (\ref{quantity}),
\noindent
we must have $V \le \hat{V}$. \\ 

\indent We now consider the reverse inequality $V \ge {\hat{V}}$. By property (f), equality holds in the last argument for the ``right $\tau$.'' Although this ``right $\tau$'' does not always exist in such problems, it does for our problem; the ``right $\tau$'' is the first entry time of the underlying Markov process $(Z,S)$ in the set where equality holds in (e). Further, this ``right $\tau$'' is a particular stopping time that is ``approximable by uniformly bounded ones.'' Larry Shepp used the phrase ``approximable by uniformly bounded ones'' to denote that we can take the  ``right $\tau$'' at which the equality is attained, approximate this ``right $\tau$'' with  ``right $\tau$'' $\wedge\, \,n$ for $n \ge 1$, and then proceed to pass to the limit for $n$. This is valid in our setting since the ``right $\tau$'' has finite expectation. When property (f) holds, and when equality holds in (d), $Y$ will be a local martingale up to the first entry of the underlying Markov process $(Z,S)$ into the the set where equality holds in (e).  Since the ``right $\tau$'' has a finite expectation, we may invoke the standard form of Doob's stopping theorem for bounded stopping times, as in \cite{Doob}. Thus,

\bneqn \label{coolone}
\expesub{\{x,r,s_1,\ldots,s_n\}}{S_1(\tau)+\ldots+S_n(\tau)- C \tau} \geq {\hat{V}}\parens{x,r,s_1,\ldots,s_n,C}.
\eneqn
and one can optimize over $\tau$ on both sides. The reverse inequality $V \ge {\hat{V}}$ thus holds and thus
$V \equiv {\hat{V}}$, completing the proof.

\end{proof}
If we can find the right ${\hat{V}}$ satisfying properties (a)-(f), we then know that

\bneqn \label{thetaeq}
A_n(C) := V\parens{O,r,0,\ldots,0;C} = \frac{\theta_n}{C} \label{eqn:first},
\eneqn

\noindent
where $\theta_n$ is a number independent of $C$. $V$ must be of the form
$\frac{\theta_n}{C}$ because a scaling argument allows us to reduce the
problem to any one value of $C$. This is because we will show that

\bneqn
V\parens{x,r,s_1,\ldots,s_n;C} = \frac{1}{C}V\parens{Cx,r,Cs_1,\ldots,Cs_n;1}. \label{eqn:C}
\eneqn
Note that if we start at $x = O$ and $s_1 = \ldots = 0$ then above form for
$A_n(C)$ is obtained.
Let 
\beqn
S(\tau)\triangleq S_1(\tau)+...+S_n(\tau).
\eeqn For any $C$ and any $\tau$, $\expe{S(\tau)} \le A_n(C)+C\expe{\tau}$. If we specify $m = \expe{\tau}$ for any fixed stopping time $\tau$, then
we will obtain the best upper bound by minimizing over $C$, which is

\beqn
\expe{S(\tau)} \le \inf\parens{\frac{\theta}{C} + C m}.
\eeqn
The infimum is attained at $C = \sqrt{\frac{\theta_n}{m}}$ and gives the bound
$C_n = 2\sqrt{\theta_n}$ for any $n$. Thus we need only find $V(O;C)$ for
any one value of $C$. 

\subsection{Solution for n=0,1,2}

\begin{corollary}
$C_0=1$.
\end{corollary}

\begin{proof}
For $n=0$, consider the function
\beqn
{\hat{V}}\parens{x,s,C} = C\parens{\parens{x-s+\frac{1}{2C}}^+}^2 + s.
\eeqn

\noindent
We note that properties (a)-(f) hold, and so for $x = s = 0$, and for any $C > 0$

\bneqn
\expe{S_\tau} \le C\expe{\tau} + \frac{1}{4C}.
\eneqn

\noindent
Minimizing over $C$, i.e., taking $C = \frac{1}{2\sqrt{\expe{\tau}} }$, as above
for any $\tau$, we obtain the inequality 

\bneqn
\expe{S(\tau)} \le \sqrt{\expe{\tau}}
\eneqn
for all
$\tau$, i.e., $C_0 = 1$.
\end{proof}

\begin{corollary}
$C_1=\sqrt{2}$.
\end{corollary}

\begin{proof}
For $n=1$, the right ${\hat{V}}$ is given by:

\bneqn
{\hat{V}}\parens{x,s,C} = Cx^2+ \frac{1}{2C}, \quad 0 \le x \le s \le \frac{1}{2C};
\eneqn

\beqn
{\hat{V}}\parens{x,s,C} = C\parens{\parens{x-s+\frac{1}{2C}}^+}^2 +s, \quad 0 \le x \le s, \,\,s > \frac{1}{2C}.
\eeqn
We use the above argument to see that $A_2(C) = \frac{1}{2C}$
and $\theta_1 = \frac{1}{2}$ and so $C_1 = \sqrt{2}$.
\end{proof}

\begin{corollary}
$C_2=\sqrt{3}$.
\end{corollary}

\begin{proof}
For $n=2$, ${\hat{V}}$ is, for $i \ne j$ \textit{and} $s_1+s_2 \le \frac{1}{C}$,
\bneqn \label{difficulteqn}
{\hat{V}}(x,r,s_1,s_2,C) = Cx^2 - Cx(s_i-s_j) + \frac{C(s_1^2+s_2^2)}{2} + \frac{3}{4C}, \quad 0 \le x \le s_i.
\eneqn
We further simplify equation (\ref{difficulteqn}) as follows:

\small
\bneqn
{\hat{V}}(x,r,s_1,s_2,C) = C\parens{\parens{x -\parens{s_i-\frac{1}{2C}}}^+}^2 + C\parens{\parens{-x -\parens{s_j-\frac{1}{2C}}}^+}^2 + s_1+s_2,
\eneqn
\normalsize
where $0 \le x \le s_i, s_1+s_2 \ge \frac{1}{C}$. We can use the above argument to see that with $V(O;C) = \frac{3}{4C}$
we arrive at $C_2 = \sqrt{3}$.
\end{proof}

\section{$n=3$ and beyond}\label{sec:formal_statement_results}

At present, we possess a non-trivial but ultimately incomplete strategy for addressing the case $n=3$. Our strategy is to develop the ``correct'' nonlinear Fredholm equation in order that we may reduce the problem to that of a nonlinear integral recurrence. Based on simulation approaches, we \textit{conjecture} the following about the constant:

\begin{conjecture}
The $\sqrt{n+1}$ pattern for the spider constant does not hold for $n=3$.
\end{conjecture}
\noindent Further, it is likely that the spider constant for $n=3$ is not an elementary number.

\section{Final Remarks}We are hopeful of a solution to the general $n$ case for the Dubins spider and maintain that our proposed dynamic programming approach constitutes the most tractable direction for solving the problem, for the following reasons:  1) The use of linear programming would be infeasible because the approximate linear programming would be large and unwieldy, making accurate numerics impossible. 2) Bellman's dynamic programming method seems intractable for the same reason as that of using linear programming. 3) The more standard method of dynamic programming, namely that of guessing a candidate for an optimal strategy, calculating the reward function of the strategy, and proving its excessiveness (as most recently done by \cite{Dubins2}) was unsuccessful in obtaining the general solution.

\section*{Acknowledgments}
First and foremost, I am indebted to my mentor, Professor Larry Shepp, for his extraordinary support, for introducing me to this literature, and for his enormously insightful conversations about this problem. I am also indebted to my colleague Professor Goran Peskir for his excellent inspiration and insight, particularly regarding the proof of Theorem 3.1.  I am grateful to Quan Zhou for his invaluable help in producing the figure in this note as well as for his careful reading of the manuscript. I thank Professor David Gilat and Professor Isaac Meilijson for their detailed input. I thank Professor Ton Dieker for his helpful comments. Finally, I am tremendously grateful to an anonymous referee whose very helpful comments enormously improved the quality of this work.

\bibliographystyle{plain}
\bibliography{VF}

\end{document}